\documentclass[11pt]{article}
\usepackage{geometry}                
\usepackage{graphicx}
\usepackage{amssymb}
\usepackage{epstopdf}
\usepackage{amsmath,amsthm,amscd,amssymb}
\usepackage{latexsym}
\numberwithin{equation}{section}

\theoremstyle{plain}
\newtheorem{theorem}{Theorem}[section]
\newtheorem{lemma}[theorem]{Lemma}

\newtheorem{conjecture}[theorem]{Conjecture}

\theoremstyle{definition}

\theoremstyle{remark}

\newtheorem{case[theorem]}{Case}

\newcommand{\ZZZ}{\mathcal{Z}}

\title{On additive doubling and energy}
\author{Nets Hawk Katz and Paul Koester}

\begin{document}

\maketitle

\begin{abstract} We discuss some ideas related to the polynomial Freiman-Ruzsa conjecture. We show
that there is a universal $\epsilon>0$ so that any subset of an abelian group with
subtractive doubling $K$ must be polynomially related to a set with additive energy at least
$\frac{1}{K^{1-\epsilon}}$. This means that the main difficulty in proving the polynomial
Freiman-Ruzsa conjecture consists in studying sets whose energy is greater than that implied
by their doubling. One example is a generalized arithmetic progression of high dimension which
cannot occur in finite characteristic.
\end{abstract}

\tableofcontents

\section{Introduction}

We are interested in studying the structure of a subset $A$ of some
additive group, which satisfies a doubling condition
$$|A+A| \leq K|A|.$$
One type of result along these lines is the class of so-called ``Freiman
theorems".

\begin{theorem} \label{FrFin} There is a function $f$ from the positive reals to the positive
reals so that if
 $A \subset F_p^n$, where $p$ is a prime and $F_p$ the field
of $p$ elements and
$$|A+A| \leq K|A|,$$
then there exists a subspace $H$ of $F_p^n$ so that $A$ is contained
in some translate of $H$ and
$$|H| \leq f(K) |A|.$$
\end{theorem}

\vskip.125in

Various authors have proved results along these lines, \cite{F73}, \cite{R99}, \cite{GT07a}, \cite{GT07b}. There is no
better possible $f$ than $p^K$ and this was essentially obtained by
Green and Tao \cite{GT07b} in the setting where $p=2$.

One may hope for a better result, namely the polynomial
Freiman-Ruzsa conjecture:

\begin{conjecture} \label{pfrfin} There is a universal constant $C>0$, so that if $A \subset F_p^n$
with
$$|A+A| \leq K|A|,$$
then there exists $z \in F_p^n$ and $H$ a subspace so that
$$|H| \leq K^C |A|,$$
and
$$|A \cap (z+H)| \geq K^{-C} |A|.$$
\end{conjecture}

The best results towards this conjecture in the case $p=2$ were obtained by Green and Tao\cite{GT07a} using
an energy-incrementation method. One obtained a subset of $A$ well situated with respect to some
subspace $H$ by taking subsets of $A$ with gradually larger additive energy. The result they obtained
was a subset $H$ with
$$|A \cap (z+H)| \geq e^{-C \sqrt{K}} |A|,$$
This was obtained by showing that given a set in $F_2^n$ with additive energy $E(A) > {1 \over K}$, one can
find a hyperplane  $P$ so that the additive energy $E(A \cap P)$ was at least ${1 \over K} + {1 \over K^{{3 \over 2}}}$.

This argument of Green and Tao at first seemed to us ripe for improvement since it
is rather inefficient for the best examples that we knew. In characteristic 2, the best example we know for a set
with additive energy ${1 \over K}$ is a set $A$ with additive doubling $K$ obtained from a subspace $H$ and a random set
$R$ of cardinality $K$ by
$$A=R+H.$$
In this case, if we take a hyperplane $P$ which contains $H$ we would expect to catch only half of the
random set and double the energy straight away. The difficulty which is at the heart of the polynomial
Freiman Ruzsa conjecture is to decompose a set $A$ with small doubling into its structured and
random parts.

In the case that $A=R+H$, we observed that one could identify $H$ by fixing an element $t$ of $A+A$,
defining $A[t]$ as that subset of $A$ which lies in $t+A$ and then calculating $A[t]+A[t]$. We wanted to
see how much of this structure remains for a general set with small doubling constant.

In the end, we failed to prove the polynomial Freiman Ruzsa conjecture, but we showed the following.
We supposed that $|A+A| \leq K |A|$ and we discovered that there is a large subset of $A+A$
which has additive energy at least $K^{- {36 \over 37}}.$ This could not be used in an induction
because additive energy is weaker than small sum set. However, we find it an interesting and
fundamental result in its own right and it is the topic of this paper.

We remark that the failure of the induction to close is a fundamental failure of our method of proof, which
works equally well independent of the characteristic. In characteristic zero, there is a basic example of
a set with large additive energy whose energy is significantly greater than the reciprocal of its doubling
constant, namely an arithmetic progression with large dimension. This example needs to appear in any
statement of polynomial Freiman Ruzsa for characteristic zero. However, we don't believe that anything
like it can appear in the finite characteristic setting. If we could show that the extreme examples
for polynomial Freiman Ruzsa in finite characteristic were not of this nature, but rather
had subtractive doubling comparable to the reciprocal of the additive energy then our result
would yield the conjecture.

\vskip.125in

\section{Preliminary Lemmas}
In the above discussion we restricted attention to the $p=2$ case but our method applies equally well when $A$ is
contained in an arbitrary finite abelian group, provided difference sets are used in place of sumsets.
We give our argument in this general setting.

\vskip.125in

We first give two useful formulations of the Cauchy-Schwarz
inequality.

\begin{lemma} \label{CS1} Let $A$ and $B$ be sets and let $f:A \longrightarrow B$. Then
$$|\{ (a,a^{\prime}): f(a)=f(a^{\prime})\}| \geq {\frac{|A|^2}{|B|}}.$$
\end{lemma}

\vskip.125in

\noindent{\it Proof.}
We calculate that
$$|\{ (a,a^{\prime}): f(a)=f(a^{\prime})\}| =\sum_{b \in B} |f^{-1}(b)|^2,$$
and applying the Cauchy-Schwarz inequality, we find
$$ \sum_{b \in B} |f^{-1}(b)|^2 \geq \frac{1}{|B|} (\sum_{b \in B} |f^{-1}(b)|)^2 =\frac{|A|^2}{|B|}.$$
\qed

\vskip.125in

\begin{lemma} \label{CS2} Let $B_1$ and $B_2$ be sets. For each $\alpha \in B_1$, let
$B_{\alpha} \subset B_2$ with $|B_{\alpha}| \geq \rho |B_2|$.
Then
$$\sum_{\alpha \in B_1} \sum_{\beta \in B_1} |B_{\alpha} \cap B_{\beta}| \geq \rho^2 |B_1|^2 |B_2|.$$
\end{lemma}

\vskip.125in

\noindent{\it Proof.} By assumption
$$\sum_{x \in B_2} \sum_{\alpha \in B_1} 1_{B_{\alpha}}(x) \geq \rho |B_1| |B_2|.$$
Applying Cauchy-Schwarz, we find
$$\sum_{x \in B_2} (\sum_{\alpha \in B_1} 1_{B_{\alpha}}(x))^2 \geq \rho^2 |B_1|^2 |B_2|.$$
But this is precisely what we were to prove.
\qed

\vskip.125in

We will make frequent use of the dyadic pigeonhole principle
throughout.  We give two formulations of this principle here.

\vskip.125in

\begin{lemma}\label{DP1}  Let $f$ be a positive real valued function defined on a finite set $S$ obeying
$$\theta\|f\|_{\infty} \leq f(s) \leq \|f\|_{\infty}$$ for all $s\in S.$  There exists
$0\leq j \leq \log_2{\theta^{-1}}$ so that
$$|\{s\in S: 2^{-j-1}\|f\|_{\infty} < f(s) \leq 2^{-j}\|f\|_{\infty}\}|
\geq \frac{1}{1 - \log_2{\theta}} |S|$$
\end{lemma}

\noindent{\it Proof.}  Let $k\in \mathbb{N}$ be the largest positive
integer so that $\theta \leq 2^{-k+1}.$  For each $0\leq j \leq
k-1,$ define
$$S_j = \{s\in S: 2^{-j-1}\|f\|_{\infty} < f(s)\leq 2^{-j}\|f\|_{\infty}\}.$$
Since $$S = \bigcup_{j=0}^{k-1} S_j,$$  there must exist some $0\leq
j \leq k-1$ so that $|S_j| \geq \frac{1}{k}|S|.$  The result follows
since\\ $k \leq 1 + \log_2{\frac{1}{\theta}}.$ \qed

\vskip.125in

\begin{lemma}\label{DP2} Let $f$ be a nonnegative valued function on a finite set $S,$ not identically zero,
and define $\theta\in (0,1]$ by
$$\frac{1}{|S|}\sum_{s\in S} f(s) = \theta \|f\|_{\infty}.$$  Then there exists a nonnegative integer
$0\leq k \leq \log_{2}{\frac{2}{\theta}}$ so that
$$|\{s \in S : 2^{-k-1}\|f\|_{\infty} < f(s) \leq 2^{-k}\|f\|_{\infty}\}| \geq
\frac{2^{k-1}\theta}{\log_{2}{(\frac{4}{\theta})}}|S|.$$
\end{lemma}

\noindent{\it Proof.} Let $S^\prime = \{s\in S: f(s) \geq
\frac{\theta}{2}\|f\|_{\infty}\}.$  Since
$$\sum_{s\in S\setminus S^{\prime}}f(s) < \frac{\theta}{2}
\|f\|_{\infty}|S|,$$
then
$$\frac{\theta}{2}\|f\|_{\infty}|S| \leq \sum_{s\in S^{\prime}}f(s) \leq \|f\|_{\infty}|S^{\prime}|,$$
and therefore $|S^{\prime}|\geq \frac{\theta}{2} |S|.$

Then
$$\frac{\theta}{2}|S| \|f\|_{\infty} \leq \sum_{j=0}^{1 + \log_{2}{\theta^{-1}}}\sum_{s\in S_j} f(s)$$
and so there exists $0 \leq k \leq \log_{2}{\frac{2}{\theta}}$ so
that
$$\frac{\theta}{2\log_{2}\frac{4}{\theta}}|S| \|f\|_{\infty} \leq
\sum_{s\in S_k} f(s) \leq 2^{-k}|S_k|\|f\|_{\infty}$$

\qed

\vskip.125in

Let $\ZZZ$ be an abelian group and $A\subset \ZZZ$ a finite subset.  The \emph{energy of $A$} is defined by
$$E(A) = |A|^{-3} \sum_{z\in \ZZZ} |(z+A)\cap A|^2 = |A|^{-3}|\{(a_1,a_2,a_3,a_4)\in A^4 : a_1 - a_2 = a_3 - a_4\}|.$$

Our main lemma shows that if one set is essentially invariant under
another set, then at least of the sets has large energy.

\begin{lemma} \label{Energy} Let $0< \rho \leq 1$ and suppose $B_1$ and $B_2$ are two subsets of $\ZZZ,$ and suppose
$$B_1 \subset \{z\in \ZZZ : |(z+B_2) \cap B_2| \geq \rho |B_2|\}.$$
Then
$$E(B_1) \geq \frac{\rho^4}{16(\log_{2}{(\frac{4}{\rho^2})})^2} \frac{|B_1|}{E(B_1)|B_2|}$$
\end{lemma}

\noindent{\it Proof.}
Applying Lemma \ref{CS2} we obtain
$$\sum_{b\in B_1} \sum_{b^{\prime}\in B_1} |(b+B_2)\cap (b^{\prime}+B_2)\cap B_2| \geq \rho^2 |B_1|^2|B_2|$$
By Lemma \ref{DP2}, we find $0\leq k \leq
\log_{2}{\frac{2}{\rho^{2}}}$ and $\Omega\subset B_1^2$ so that
$$|\Omega|\geq \frac{2^{k-1}\rho^2}{\log_{2}{(\frac{4}{\rho^2})}} |B_1|^2$$
and
$$2^{-k-1}|B_2| \leq |(b+B_2)\cap (b^{\prime}+B_2) \cap B_2| \leq 2^{-k}|B_2|$$
for all $(b,b^{\prime})\in \Omega.$

Define $S=\{z\in \ZZZ : |(z+B_2) \cap B_2| \geq 2^{-k-1} |B_2|\}.$  Then
$$|B_2|^3 E(B_2) = \sum_{z\in \ZZZ}|(z+B_2) \cap B_2|^2$$
$$\geq \sum_{z\in S}|(z+B_2) \cap B_2|^2 \geq |S|2^{-2k-2}|B_2|^2,$$
giving the bound $$|S|\leq 2^{2k+2}|B_2|E(B_2).$$

On the other hand,
$$|B_1|^3 E(B_1) = |\{(a_1,a_2,a_3,a_4)\in B_1^4: a_1-a_2 = a_3-a_4\}|$$
$$\geq |\{(g,g^{\prime})\in \Omega^2 : -(g) = -(g^{\prime})\}| \geq \frac{|\Omega|^2}{|-(\Omega)|}$$
$$\geq \frac{|\Omega|^2}{|S|} \geq \frac{2^{2k-2}\rho^4}{(\log_{2}
{(\frac{4}{\rho^2})})^2} \frac{|B_1|^4}{2^{2k+2}|B_2|E(B_2)},$$
which gives the desired lower bound on $E(B_1).$ \qed

\section{Main result}

Given a set $\Omega\subset \ZZZ^2$, we define $-:\Omega\rightarrow
\ZZZ$ by $-(a,b) = a-b.$  Given $A\subset \ZZZ$ and $t\in A-A,$ we
write
$$A[t] = (A+t)\cap A.$$
Then  $|A[t]| = |\{(a,b) \in A^2: a-b = t\}|.$

We let $K$ be fixed throughout.  For two quantities $A$ and $B,$ we
write $A\lesssim B$ if there exists a constant $C$, independent of
$K$, so that $A\leq C B.$  We write $A\sim B$ if $A\lesssim B$ and
$B\lesssim A.$  We write $A\lessapprox B$ if for each $\delta>0$
there exists a constant $C_{\delta}>0,$ independent of $K$, so that
$$A\leq C_{\delta}K^{\delta}B.$$  We write $A\approx B$ if $A\lessapprox B$ and
$B\lessapprox A.$

\begin{theorem}  There exists universal constants $\epsilon>0$ and $C>0$ so that if $A$ is a finite subset
of an abelian group $\ZZZ$ satisfying $|A - A| = K|A|,$ then there exists $A^\prime \subset A - A$ satisfying
$$ |A^{\prime}| \geq K^{-C}|A|$$
and
$$ E(A^{\prime}) \gtrapprox \frac{1}{K^{1-\epsilon}}.$$
Furthermore, we can take $\epsilon = \frac{1}{37}.$
\end{theorem}


\vskip.125in

The strategy of the proof is to investigate several reasonable
candidates for $A^{\prime}.$  We will always assume these candidate
sets do not satisfy the energy bound required in the theorem, since
we would be done otherwise.  But if a candidate set fails the
desired energy bound we then use Lemma \ref{Energy} to put a lower
bound on the energy of a related set.  Continuing in this manner, we
eventually find a set $A^{\prime}$ as desired in the theorem.

\vskip.125in

We assume $E(A-A) < \frac{1}{K^{1-2\epsilon}}$ throughout, since
otherwise we could take $A^{\prime} = A-A.$ Likewise, we assume
$E(A) < \frac{1}{K^{1-\epsilon}}.$

\vskip.125in

For each $t\in A-A,$ we have
$$ |A| \leq |A[t]-A| \leq K|A|;$$
we therefore obtain an exponent $0\leq \beta(t) \leq 1$ so that
$|A[t] - A| = K^{\beta(t)}|A|.$

\begin{lemma}\label{FirstRefinement} Let $\epsilon > 0$ and suppose $A$ is a finite subset of an abelian
group $\ZZZ$ satisfying $|A-A| = K|A|$ and $E(A) \lessapprox
\frac{1}{K^{1-\epsilon}}.$  There exists a set $T \subset A-A$ and
an exponent $0 \leq \beta \leq 1$ so that
\begin{equation}\label{E:ATBig}
|A[t]| \geq \frac{|A|}{2K}
\end{equation}
and
$$|A[t] - A| \sim K^{\beta}|A|$$
for all $t\in T,$
and in addition,
$$|T| \gtrapprox K^{-\epsilon}|A-A|$$
\end{lemma}

\noindent{\it Proof.} Letting $H = \{(a,a^{\prime})\in A^2: |A[a-b]| \geq \frac{|A|}{2K}\},$ we first observe
$$\sum_{(a,a^{\prime})\in A^2\setminus H} 1 = \sum_{\substack{t\in A-A\\ |A[t]|\geq \frac{|A|}{2K}}}
|A[t]| \leq \frac{|A|}{2K}|A-A| = \frac{1}{2}|A|^2,$$
and therefore $H\geq \frac{1}{2}|A|^2.$

On the other hand,
$$\frac{|H|^2}{|-(H)|} \leq \sum_{t\in -(H)}|A[t]|^2 \lessapprox \frac{|A|^3}{K^{1-\epsilon}},$$
and therefore
$$|-(H)|\gtrapprox \frac{K^{1-\epsilon}|H|^2}{|A|^3} \geq K^{1-\epsilon}|A| = K^{-\epsilon}|A-A|.$$
We have \eqref{E:ATBig} for all $t\in -(H).$

For each $t\in -(H),$ $|A[t] - A| \sim K^{\beta(t)}|A|;$  by Lemma
\ref{DP1}, there exists $T\subset -(H)$ and $0\leq \beta\leq 1$ so
that $|A[t] - A| \sim K^{\beta}|A|$ for all $t\in T,$ and
$|T|\gtrapprox |-(H)|.$  The theorem follows since $|-(H)|\gtrapprox
K^{-\epsilon}|A-A|.$ \qed

\vskip.125in

\subsection{The Large $\beta$ case: ($\beta \geq \frac{1}{2} +
\frac{7}{4}\epsilon$)}
%

\begin{lemma} There exists $0\leq \alpha\leq \beta$ and a set $X\subset A-A$ so that
$|X|\gtrapprox K^{-\alpha}|A-A|$ and
$$|\{t \in T: x \in A[t] - A\}| \gtrsim K^{\alpha + \beta - \epsilon}|A|$$ for each $x \in X.$
\end{lemma}

\noindent{\it Proof.}
$$\sum_{t\in T}|A[t]-A| \gtrapprox K^{\beta - \epsilon}|A||A-A|,$$
but
$$\sum_{t\in T}|A[t] - A| = \sum_{x\in A-A}|\{t \in T: x \in A[t] - A\}|.$$
The conclusion follows from Lemma \ref{DP2}.  \qed

\vskip.125in

Suppose $x\in A[t] - A$ for some $t\in T\subset A-A.$  Then
$$x \in t+A-A$$
and so
$$t\in x + A-A.$$  Thus, $t\in (x+A-A)\cap(A-A).$  In particular, if $x\in X,$ then
$$|(x + A - A)\cap (A - A)| \gtrsim K^{\alpha+\beta - \epsilon}|A| = K^{\alpha + \beta - \epsilon - 1}|A-A|,$$
since $x\in A[t] - A$ for at least $K^{\alpha + \beta - \epsilon}|A|$ values of
$t.$  Applying Lemma \ref{Energy} with $B_1 = X,$ $B_2 = A - A,$ and
$\rho = K^{\alpha + \beta - \epsilon - 1},$ we obtain
$$E(X) \gtrapprox K^{4\alpha + 4\beta - 4\epsilon - 4}K^{1-2\epsilon}K^{\alpha} \geq K^{4\beta - 6\epsilon - 3}.$$
If $\beta \geq \frac{1}{2} + \frac{7}{4}\epsilon$ then
$$E(X) \gtrapprox \frac{1}{K^{1-\epsilon}},$$
giving the conclusion of the Theorem with $A^\prime = X.$  This
proves the theorem provided $\beta \geq\frac{1}{2} +
\frac{7}{4}\epsilon.$

\vskip.125in

\subsection{The Small $\beta$ case. ($\beta < \frac{7}{12} -
\frac{4}{3}\epsilon$)}


\begin{lemma}For each $t\in T,$ there exists a set $G_1(t)\subset A[t]^2$ and $0\leq \alpha(t) \leq \beta$ so that
$$|G_1(t)| \gtrsim K^{-\alpha(t)}|A[t]|^2$$
and
$$|(a-a^\prime + A)\cap A|\gtrapprox K^{\alpha(t) - \beta}|A|$$
for all $(a,a^{\prime}) \in G_1(t).$
\end{lemma}
\noindent{\it Proof.}
We begin by observing
$$\sum_{(a^\prime,a)\in A[t]^2} |(a^{\prime} - A)\cap(a - A)| \geq K^{2\beta}|A[t]|^2|A[t] - A| = K^{\beta}|A[t]|^2|A|.$$
By Lemma \ref{DP2}, we obtain $0\leq \alpha(t) \leq \beta$ and a set
$G_1(t)\subset A[t]^2$ so that
$$|G_1(t)|\gtrapprox K^{-\alpha(t)}|A[t]|^2$$
and
$$|(a^{\prime} - A)\cap (a-A)| \gtrapprox K^{\alpha(t) + \beta}|A|.$$
for all $(a,a^\prime) \in G_1(t).$  But
$$|(a^{\prime} - A)\cap (a-A)| = |(A - a^{\prime})\cap (A - a)| = |(a-a^{\prime} + A)\cap A|.$$
\qed

\vskip.125in

By Lemma \ref{DP1}, we may find $0\leq \alpha \leq \beta$ and
$T^{\prime}\subset T$ so that $|T^{\prime}|\gtrapprox |T|$ and so
that $|\alpha(t) - \alpha| \leq \log_K{2}$ for all $t\in
T^{\prime}.$

\vskip.125in

\begin{lemma} For each $t\in T^{\prime},$  there exists $\gamma(t)\geq 0 $ and $G(t)\subset G_1(t)$ so that
$$|G(t)|\gtrapprox |G_1(t)|,$$
$$|-(G(t))| \gtrapprox K^{\gamma(t) - \alpha}|A[t]|$$
$$|(A[t])[x]| \sim K^{-\gamma(t)}|A[t]|$$
for all $x\in -(G(t))$, and
$$E(A[t])\gtrapprox K^{-\alpha + \gamma(t)}$$
\end{lemma}

\noindent{\it Proof.}
Let $G^{\prime}(t) = \{(a,b)\in G_1(t) : |(A[t])[a-b]|\geq \frac{|G_1(t)|}{2|-(G_1(t))|}\}.$  Then
$$\sum_{(a,b)\in G_1(t)\setminus G^{\prime}(t)}1 = \sum_{x\in -(G_1(t)\setminus G^{\prime}(t))}
|(A[t])[x]| \leq \frac{|G_1(t)|}{2},$$
and so $|G^{\prime}(t)| \geq \frac{|G_1(t)|}{2}.$

Next, we observe
$$\frac{|G_1(t)|}{2|-(G_1(t))|} \geq \frac{|G_1(t)|}{2|A[t] - A|} \gtrapprox \frac{K^{\alpha}
|A[t]|^2}{K^{\beta}|A|} = K^{-1-\alpha - \beta}|A[t]|.$$
Coupled with $|(A[t])[x]| \leq |A[t]|,$ we conclude there exists $G(t)\subset G^{\prime}(t)$
and $\gamma(t)\geq 0$ so that $|G(t)|\gtrapprox |G^{\prime}(t)|$ and $|(A[t])[x]|\sim K^{\gamma(t)}|A[t]|$
for all $x\in -(G(t)).$

Next, we observe
$$|G(t)| = \sum_{x\in -(G(t))}|(A[t])[x]| \approx K^{-\gamma(t)}|A[t]||-(G(t))|,$$
whereas we also see
$$|G(t)|\gtrapprox |G_1(t)| \gtrapprox K^{-\alpha}|A[t]|^2,$$
giving the desired bound on $|-(G(t))|.$

Finally,
$$|A[t]|^3E(A[t]) = \sum_{x\in A[t] - A[t]}|(A[t])[x]|^2 \geq \sum_{x\in -(G(t))}|(A[t])[x]|^2
\gtrapprox K^{\gamma(t) - \alpha}|A[t]|K^{-2\gamma(t)}|A[t]|^2,$$
giving the desired bound on $E(A[t]).$ \qed

\vskip.125in

Let $G(t)$ be the set found in the previous lemma.  We may assume
the stronger bound $|(-G(t))|\gtrapprox K^{1-2\alpha -
\epsilon}|A[t]|.$  Indeed, if $\gamma(t)\leq 1 - \alpha - \epsilon$
then
$$E(A[t]) \gtrapprox \frac{1}{K^{1-\epsilon}};$$
we thus obtain the conclusion of the theorem taking $A^{\prime} =
A[t]$.  We may therefore assume $\gamma(t) > 1 - \alpha - \epsilon,$
hence
$$|(-G(t))|\gtrapprox K^{1 - 2\alpha - \epsilon}|A[t]| \gtrapprox K^{-1-2\alpha - \epsilon}|A-A|$$ for
all $t\in T^{\prime}.$

\vskip.125in

Now define
$$X^{\prime} = \bigcup_{t\in T^{\prime}}-(G(t))$$
and
$$g(x) = |\{t\in T^{\prime} : x\in -(G(t))\}|.$$

\begin{lemma}  There exists $X\subset X^{\prime}$ and $\eta\geq 0$ so that
$$|X| \gtrapprox K^{-1-2\alpha-2\epsilon + \eta}|A-A|$$
and
$$ g(x) \sim K^{-\eta}|A-A|$$
for all $x\in X.$
\end{lemma}
\noindent{\it Proof.} We first observe $|X^{\prime}|\leq |A-A|.$
Then
$$\frac{1}{|X^{\prime}|}\sum_{x\in X^{\prime}} g(x) = \frac{1}{|X^{\prime}|}\sum_{t\in T^{\prime}} |-(G(t))|
\gtrapprox
\frac{|T^{\prime}|}{|X^{\prime}|}K^{-1-2\alpha-\epsilon}|A-A|$$
$$ \gtrapprox \frac{K^{-\epsilon}|A-A|}{|A-A|}K^{-1-2\alpha-\epsilon}|A-A|
\gtrapprox K^{-1-2\alpha - 2\epsilon}|A-A|.$$
Also, $g(x)\leq |T^{\prime}| \lessapprox K^{-\epsilon}|A-A| \leq
|A-A|.$ The result now follows from Lemma \ref{DP2}. \qed

\vskip.125in

If $x\in X$ then $g(x)\gtrsim K^{-\eta}|A-A|,$ so $x\in -(G(t))$ for
at least $\sim K^{-\eta}|A-A|$ values $t\in T^{\prime}.$ Therefore,
$x\in A[t] - A \subset t + A-A$ for at least $\sim K^{-\eta}|A-A|$
values of $t\in T^{\prime},$  so $t\in x + A - A$ for each of these
$t.$ But each such $t$ is also in $A-A,$ so
$$|(x + A-A)\cap (A-A)| \gtrapprox K^{-\eta}|A-A|$$ for each $x\in X.$

Applying Lemma \ref{Energy} with $B_1 = X, B_2 = A-A$ and $\rho =
K^{-\eta},$ we obtain
$$E(X) \gtrapprox K^{-4\eta} K^{1-\epsilon} K^{-1-2\alpha-2\epsilon + \eta} = K^{-2\alpha - 3\epsilon -3\eta}.$$
If $E(X)\gtrapprox \frac{1}{K^{1-\epsilon}},$ we obtain the
conclusion of the theorem by taking $A^\prime = X;$  we may
therefore assume
$$\eta > \frac{1}{3}(1 - 2\alpha - 4\epsilon).$$
Our lower bound on $X$ is now
$$|X| \gtrapprox K^{\frac{1}{3}-\frac{8}{3}\alpha-\frac{10}{3}\epsilon}|A|.$$
But for each $t\in T^{\prime}$ we have
$$|(a-a^\prime + A)\cap A|\gtrapprox K^{\alpha + \beta}|A|$$
for all $(a,a^{\prime})\in G(t).$  Therefore,
$$|(x + A)\cap A|\gtrapprox K^{\alpha - \beta}|A|$$
for all $x\in X.$

Using Lemma \ref{Energy} with $B_1 = X, B_2 = A$ and $\rho =
K^{\alpha - \beta}$ we obtain
$$E(X) \gtrapprox K^{4(\alpha - \beta)}K^{1-\epsilon}K^{\frac{1}{3}-\frac{8}{3}\alpha-\frac{10}{3}\epsilon}
= K^{\frac{4}{3} + \frac{4}{3}\alpha - 4\beta - \frac{13}{3}\epsilon}
\geq K^{\frac{4}{3} - 4\beta - \frac{13}{3}\epsilon}.$$
If $E(X)\gtrapprox \frac{1}{K^{1-\epsilon}}$ we obtain the
conclusion of the theorem upon taking $A^{\prime} = X.$  We may
therefore assume
$$\frac{4}{3} - 4\beta - \frac{13}{3}\epsilon < -1 + \epsilon, $$
so
$$\frac{7}{12} - \frac{4}{3}\epsilon < \beta.$$
But we already proved the theorem in the case $\beta \geq
\frac{1}{2} + \frac{7}{4}\epsilon.$ Taking $\epsilon =
\frac{1}{37},$ the $\beta
> \frac{7}{12} - \frac{4}{3}\epsilon$ range is now subsumed by the
$\beta \geq \frac{1}{2} + \frac{7}{4}\epsilon$ case.  This concludes
the proof of our main theorem.  \qed

\end{document}